\documentclass[11pt]{amsart}
\usepackage{graphicx}
\usepackage[usenames]{color}
\usepackage{amssymb,amscd}
\newcommand\Z{\mathbb Z}
\newcommand\Q{\mathbb Q}
\newcommand\R{\mathbb R}

\newcommand\Rinf{R_{\infty}}

\newcommand\vp{\varphi}

\newcommand\ph{\varphi}
\newtheorem{theorem}{Theorem}[section]

\newtheorem{lemma}[theorem]{Lemma}

\newtheorem{definition}[theorem]{Definition}

\title{A note on twisted conjugacy and generalized Baumslag-Solitar groups}
\author{Jennifer Taback}
\address{Department of Mathematics,
Bowdoin College, Brunswick, ME 04011} \email{jtaback@bowdoin.edu}
\author{Peter Wong}
\address{Department of Mathematics, Bates College, Lewiston, ME 04240} \email{pwong@bates.edu}
\thanks{The first author acknowledges support from
NSF grant DMS-0604645, and would like to thank Kevin Whyte for many
useful conversations about this paper.  Both authors would like to
thank Gilbert Levitt for bringing reference \cite{L} to their
attention and Christopher Cashen for helpful comments regarding the proof of Theorem 3.1.} 

\begin{document}

\maketitle

\begin{abstract}
A generalized Baumslag-Solitar group is the fundamental group of a
graph of groups  all of whose vertex and edge groups are infinite
cyclic.  Levitt proves that any generalized Baumslag-Solitar group
has {\em property $\Rinf$}, that is, any automorphism has an
infinite number of twisted conjugacy classes.  We show that any
group quasi-isometric to a generalized Baumslag-Solitar group also
has property $\Rinf$.  This extends work of the authors proving that
any group quasi-isometric to a solvable Baumslag-Solitar $BS(1,n)$
group has property $\Rinf$, and relies on the classification of
generalized Baumslag-Solitar groups given by Whyte.
\end{abstract}

\section{Introduction}
We say that a group $G$ has {\em property $\Rinf$} if any
automorphism  $\vp$ of $G$ has an infinite number of twisted
conjugacy classes.  Two elements $g_1,g_2 \in G$ are $\vp$-twisted
conjugate if there is an $h \in G$ so that $hg_1\vp(h)^{-1} = g_2$.
The study of the finiteness of the number of twisted conjugacy
classes arises in Nielsen fixed point theory. For example, for each
$n\ge 5$, there is a compact $n$-dimensional nilmanifold $M^n$ whose
fundamental group has property $R_{\infty}$ \cite{GW3}. As a
consequence, {\it every} homeomorphism of such a manifold $M^n$ is
isotopic to a fixed point free homeomorphism. For more details and
background on fixed point theory, see \cite{B} or \cite{J}.

Recently, several authors have studied property $\Rinf$ from a
geometric perspective, where the word geometric has a variety of
interpretations.  It is proven in both \cite{FG2,TWh} that a group
which has a non-elementary action by isometries on a Gromov
hyperbolic space has property $\Rinf$, where the action is
fundamental to understanding the twisted conjugacy classes.
Recently, the authors have given a proof that the lamplighter groups
$L_n = \Z_n \wr \Z$ have property $\Rinf$ iff $(n,6) = 1$,
originally proven in \cite{GW1}, using mainly the geometry of the
Cayley graph of these groups.  Namely, the geometry of the
Diestel-Leader graph $DL(n,n)$ combined with recent results of
Eskin, Fisher and Whyte \cite{EFW1,EFW2} provides a geometric
interpretation for the twisted conjugacy classes. \cite{TW2}

A natural question to ask is whether property $\Rinf$ is {\em
geometric}, that is, invariant under quasi-isometry.  It is shown in
\cite{FG} that the Baumslag-Solitar groups $BS(m,n)$ (excepting
$BS(1,1)$) have property $\Rinf$, and in  \cite{TW1} that any group
quasi-isometric to $BS(1,n)$ also has the property.  The analogous
results are shown in \cite{TW1} for the solvable generalization
$\Gamma$ of $BS(1,n)$ given by the short exact sequence
$$1 \rightarrow \Z[\frac{1}{n}] \rightarrow \Gamma \rightarrow \Z^k \rightarrow 1$$
 and any group quasi-isometric to $\Gamma$.  However, property
 $\Rinf$ is not in general a quasi-isometry invariant.   Let $A, \ B \in
GL(2,\Z)$ be matrices whose traces have absolute value at least two.
Then $\Z^2 \rtimes_A \Z$ and $\Z^2 \rtimes_B \Z$ are always
quasi-isometric, as they are both cocompact lattices in Sol, but may
not both have property $\Rinf$ \cite{GW2}.

In this note we prove the following theorem about groups
quasi-isometric to generalized Baumslag-Solitar groups, extending
the result of Levitt \cite{L} as well as the results of \cite{TW1}.
A generalized Baumslag-Solitar group is the fundamental group of a
graph of groups all of whose vertex and edge groups are infinite
cyclic.

\noindent {\bf Theorem 3.1.} {\em Let $G$ be a finitely generated
group quasi-isometric to a non-elementary  generalized
Baumslag-Solitar group.  Then $G$ has property $\Rinf$.}

Our proofs rely on a result of Whyte \cite{W} stating that any group
$G$  quasi-isometric to a generalized Baumslag-Solitar group has one
of three forms: either $G$ is $BS(1,n)$ and the result is proven in
\cite{TW1}, $G$ is virtually $F_n \times \Z$ or $G$ is the
fundamental group of a graph of groups all of whose vertex and edge
groups are {\em virtually} infinite cyclic.

In the second case above, we rely on work of Sela \cite{S} at a
crucial step to guarantee that the quotient group we are
considering, which is a non-elementary Gromov hyperbolic group, is
Hopfian.  This quotient group is obtained by considering a {\em
quasi-action}, as an action may not exist, of the group on the
product of a tree with the real line.

The third case splits into two subcases depending on whether $G$ is
unimodular. If it is not, then Levitt's proof that any generalized
Baumslag-Solitar group has property $\Rinf$ applies verbatim to $G$.
When a generalized Baumslag-Solitar group is unimodular, Levitt
studies central elements which are necessarily elliptic, and
concludes that the center of the group is either trivial or infinite
cyclic. When $G$ is quasi-isometric to a generalized
Baumslag-Solitar group, we can only conclude that the center is {\em
virtually} $\Z$ or else the normal closure of the torsion elements
of $G$.  We then show that this cannot occur, that is, in this case
$G$ cannot be unimodular.

Thus the complete proof of Theorem \ref{thm:gbs} combines a number
of existing techniques and theorems due to Levitt \cite{L}, Kleiner
and Leeb \cite{KL}, Whyte \cite{W} and the authors \cite{TW1}, and
adds a new and interesting class to the list of groups for which
property $\Rinf$ is invariant under quasi-isometry. We note that
part of the main theorem could also be proven using a result of
Rieffel \cite{R}.

\section{Background on twisted conjugacy and quasi-isometries }
\subsection{Twisted conjugacy}
Let $\ph:\pi \to \pi$ be a group endomorphism.  We consider the
action  of $\pi$ on $\pi$ given by $\sigma \cdot \alpha \mapsto
\sigma \alpha \ph(\sigma)^{-1}$ for $\sigma, \ \alpha \in \pi$. The
orbits of this action are the {\it Reidemeister classes} of $\ph$ or
the $\ph$-twisted conjugacy classes. Denote by $R(\ph)$ the
cardinality of the set $\mathcal R(\ph)$ of $\ph$-twisted conjugacy
classes.  This number $R(\ph)$ is called the {\em Reidemeister
number} of $\ph$. When $\ph$ is the identity, $\mathcal R(\ph)$ is
the set of conjugacy classes of elements of $\pi$, and $R(\ph)$ is
simply the number of conjugacy classes.

We say that a group $G$ has {\em property $\Rinf$} if, for any  $\vp
\in Aut(G)$, we have $R( \vp) = \infty$.  The main technique we use
for computing $R(\vp)$ is as follows.  We consider groups which can
be expressed as group extensions, for example $1 \rightarrow A
\rightarrow B \rightarrow C \rightarrow 1$.  Suppose that an
automorphism $\vp \in Aut(B)$ induces the following commutative
diagram, where the vertical arrows are group homomorphisms, that is,
$\ph|A = \ph'$ and ${\overline \ph}$ is the quotient map induced by
$\ph$ on $C$:

\begin{equation}\label{short-exact}
\begin{CD}
    1   @>>> A    @>{i}>>  B @>{p}>>      C @>>> 1 \\
    @.  @V{\ph'}VV      @V{\ph}VV   @V{\overline \ph}VV @.\\
    1   @>>> A    @>{i}>>  B @>{p}>>      C @>>> 1
 \end{CD}
\end{equation}

Then we obtain a short exact sequence of sets and corresponding functions ${\hat i}$ and ${\hat p}$:
$$
\mathcal R(\ph') \stackrel{\hat i}{\to} \mathcal R(\ph) \stackrel{\hat p}{\to} \mathcal R(\overline \ph) \to 1
$$
where if $\bar 1$ is the identity element in $C$, we have
$\hat{i}(\mathcal R(\ph'))={\hat p}^{-1}([\bar 1])$,  and $\hat p$
is onto. To ensure that both $\varphi'$ and $\overline \ph$ are both
automorphisms, we need the following lemma.  Recall that a group is
{\em Hopfian} if every epimorphism is an automorphism.  The
following lemma is proven in \cite{TW2}.

\begin{lemma}\label{lemma:hopfian}
If $C$ is Hopfian, then $\ph'\in Aut(A)$ and $\overline \ph \in
Aut(C)$.
\end{lemma}

The following result is straightforward and follows from more
general results discussed in \cite{Wo}.

\begin{lemma}\label{lemma:reid}
Given the commutative diagram labeled \eqref{short-exact} above,
\begin{enumerate}
\item if $R(\overline \ph)=\infty$ then $R(\ph)=\infty$, and

\item if $C$ is finite and $R(\ph')=\infty$ then $R(\ph)=\infty$.
\end{enumerate}
\end{lemma}

\subsection{Quasi-isometries and quasi-actions}
A quasi-isometry is a map between metric spaces which distorts distance by a uniformly bounded amount, defined precisely as follows.
\begin{definition}
Let $X$ and $Y$ be metric spaces.  A map $f: X \rightarrow Y$ is a $(K,C)$-quasi-isometry, for $K \geq 1$ and $C \geq 0$ if
\begin{enumerate}
\item $\frac{1}{K} d_X(x_1,x_2) - C \leq d_Y(f(x_1),f(x_2)) \leq K d_X(x_1,x_2) +C$ for all $x_1,x_2 \in X$.
\item For some constant $C'$, we have $Nhbd_{C'}(f(X)) = Y$.
\end{enumerate}
\end{definition}

There is a notion of {\em coarse inverse} for a quasi-isometry; quasi-isometries $f$ and $g$ are coarse inverses if both compositions $f \circ g$ and $g \circ f$ are a bounded distance from the identity.  Without loss of generality we may assume that $g$ and $f$ share the same quasi-isometry constants.

The set of all self quasi-isometries of a finitely generated group $G$ is denoted $QIMap(G)$.  We form equivalence classes consisting of all quasi-isometries of $G$ which differ by a uniformly bounded amount.  This set of equivalence classes is called the {\em quasi-isometry group} of $G$ and denoted $QI(G)$.

If $G$ is a finitely generated group and $X$ is a proper geodesic metric space, we define a {\em quasi-action} of $G$ on $X$ to be a map $\Psi: G \rightarrow QIMap(X)$ satisfying the following properties for some constants $K \geq 1$ and $C \geq 0$.
\begin{enumerate}
\item For each $g \in G$, the element $\Psi(g)$ is a $(K,C)$-quasi-isometry of $X$.
\item $\Psi(Id)$ is a uniformly bounded distance from the identity, that is, $\Psi(Id)$ is the identity in $QI(X)$.
\item $\Psi(g) \Psi(h)$ is a uniformly bounded distance from $\Psi(gh)$, for all $g,h \in G$.
\end{enumerate}
It is clear that a quasi-action induces a homomorphism from $G$ into $QI(X)$.  It is important to note that all quasi-isometries $\Psi(g)$ have the same quasi-isometry constants.

When using Lemma \ref{lemma:reid} to prove that a group $G$ has property $\Rinf$, one must be able to find characteristic subgroups of $G$.  One such characteristic subgroup of $G$ is the {\em virtual center} $V(G)$.  This subgroup consists of elements of $G$ whose centralizers have finite index in $G$, that is, $V(G)=\{g\in G| [G:C(g)]<\infty\}$.  When $G$ has a quasi-action on a Cayley complex, it follows from Lemma \ref{virtual-center} below that the virtual center of $G$ is exactly the kernel of this quasi-action.  Namely, we use the fact that $g \in G$ moves all points of $G$ a uniformly bounded distance under left multiplication if and only if $g$ has finitely many conjugates, and show that the virtual center consists exactly of those group elements having finitely many conjugates.  This gives another characterization of the virtual center of $G$.  By construction, the virtual center is a characteristic subgroup of $G$.

\begin{lemma}\label{virtual-center}
Let $G$ be a group generated by a finite set $S$ which has a quasi-action on a Cayley complex $X$.  Then the virtual center of $G$ consists of those elements that move all points of $X$ a uniformly bounded distance $B$, that is, $ \ d(gx,x) \leq B$ for all $x \in X$, where distance is computed in the word metric with respect to $S$.
\end{lemma}

\begin{proof}
To prove this lemma, we use the fact that $g \in G$ moves all points of $G$ a uniformly bounded distance under left multiplication if and only if $g$ has finitely many conjugates.  This is true because there are a finite number of elements of $G$ in the ball of radius $B$ in any Cayley graph of $G$ with respect to a finite generating set.

We first show that if $g \in V(G)$, then $g$ has finitely many conjugates in $G$.
Fix an element  $g$ in the virtual center $V(G)$.  Since $[G: C(g)] < \infty$, we can write $G$ as a disjoint union of the cosets of $C = C(g)$, namely $C,a_1C,a_2C, \cdots a_nC$ for some $a_1,a_2, \cdots ,a_n \in G$.

Take any two elements from the coset $a_1C$, say $a_1c_1$ and $a_1c_2$ where $c_1,c_2 \in C$.  When we conjugate $g$ by these two elements we see that
$$ (a_1c_1)g(a_1c_1)^{-1} = a_1(c_1gc_1^{-1})a_1^{-1} = a_1ga_1^{-1},$$ where the last equality follows because $c_1 \in C(g)$, and
$$ (a_1c_2)g(a_1c_2)^{-1} = a_1(c_2gc_2^{-1})a_1^{-1} = a_1ga_1^{-1}.$$
Thus we see that the only possible conjugates of $g$ have the values $a_iga_1^{-1}$ for $i = 1, \cdots n,$ or $g$ itself.

Now assume that $g$ has finitely many conjugates in $G$.  We will show that the centralizer $C=C(g)$ has finite index in $G$.  Suppose there were infinitely many cosets of $C$, which we denote $a_1C,a_2C, \cdots$ for an infinite sequence $a_1ma_2 \cdots \in G$.  Consider the conjugates $a_i^{-1}ga_i$.  Since $g$ has finitely many conjugates, we know that infinitely many of these must be the same.  Suppose that $a_1^{-1}ga_1 = a_2^{-1}ga_2$.  This is equivalent to $(a_2^{-1}a_1)^{-1}g(a_2^{-1}a_1) = g$, so $a_2^{-1}a_1$ must be in $C(g)$, so as cosets $a_1C = a_2C$. Thus there are a finite number of cosets and $C$ has finite index in $G$.
\end{proof}

\subsection{Group actions on trees}
We assume that $G$ is a finitely generated group acting simplicially on a tree $T$, that is, this action preserves vertices and edges of $T$, without inversions.  Moreover, we require this action to be {\em minimal}, meaning that there is no proper invariant subtree under this action.

An element $g \in G$ is called {\em elliptic} if $g$ fixes a vertex in $T$, and {\em hyperbolic} otherwise.  These properties are best defined in terms of translation length, as follows.  View $T$ as a metric space by assigning each edge length one.  Define the {\em translation length} $l_g$ of an element $g \in G$ acting on $T$ to be the minimum distance between a vertex $x \in V(T)$ and its image $g \cdot x$, that is,
$$l_g = \min_{x \in V(T)} d(x,g \cdot x).$$
If $l_g = 0$ then we say that $g$ is {\em elliptic}, and $g$ has a fixed point when acting on $T$.  Otherwise $g$ is {\em hyperbolic}, and $T$ has a $g$-invariant linear subtree $A_g$, called the axis of $g$, consisting of the following set of points:
$$A_g = \{x \in T | d(x,g \cdot x) = l_g \}.$$

The notion of commensurability will also play a role in the proofs below.  There are several standard definitions of commensurability.  Two groups are abstractly commensurable if they have isomorphic finite index subgroups.  Two subgroups $H_1$ and $H_2$ of a given group $G$ are commensurable if their intersection has finite index in both subgroups.  We say that two elements $g,h \in G$ are commensurable if the subgroups they generate are commensurable in $G$.

When $G$ is a generalized Baumslag-Solitar group, all elliptic elements have infinite order and are commensurable.  For the groups we consider, in which all vertex and edge stabilizers are only virtually infinite cyclic, all infinite order elliptic elements are still commensurable, and the same is true for any finite order elliptic elements. In general, the properties of being elliptic or hyperbolic and having finite or infinite order are preserved under both conjugation and commensurability. This may not be true when considering twisted conjugacy, however.  In the case of the infinite dihedral group $D_{\infty}$, the order two elliptic elements are $\vp$-twisted conjugate to the infinite order elliptic elements where $\vp: \Z \rightarrow \Z$ is given by $\vp(x) = -x$.

\section{Twisted conjugacy and generalized Baumslag-Solitar groups}

Below we prove that any group quasi-isometric to a non-elementary
generalized  Baumslag-Solitar group has property $\Rinf$.  A
generalized Baumslag-Solitar group is a finitely generated group
which acts on a tree with all edge and vertex stabilizers infinite
cyclic, that is, the fundamental group of a graph of groups all of
whose vertex and edge groups are infinite cyclic.  A group is
non-elementary if it is not virtually cyclic.

According to \cite{W}, Theorem 0.1, any group $\Gamma$ which is a
generalized  Baumslag-Solitar group has one of three forms:
\begin{enumerate}
\item $\Gamma = BS(1,n)$ for some $n > 1$
\item $\Gamma$ is virtually $F_n \times \Z$
\item $\Gamma$ is quasi-isometric to $BS(2,3)$.
\end{enumerate}
We will consider these three cases in the proof below.  We now prove the following theorem.

\begin{theorem}\label{thm:gbs}
Let $G$ be a finitely generated group quasi-isometric to a
non-elementary  generalized Baumslag-Solitar group.  Then $G$ has
property $\Rinf$.
\end{theorem}
\begin{proof}

{\bf Case 1.} If $G$ is quasi-isometric to $\Gamma = BS(1,n)$ for
some $n > 1$ then it is proven in \cite{TW1} that $R(\vp) = \infty$
for all $\vp \in Aut(G)$.

{\bf Case 2.} If $G$ is quasi-isometric to $\Gamma$, which is
virtually $F_n \times \Z$, then $G$ itself is quasi-isometric to
$F_n \times \Z$.  The geometric model of this group, also called the
Cayley complex, which is quasi-isometric to the group, is then the
product of a tree $T$ with $\R$.  In particular, since $n \geq 2$,
the tree $T$ is not a line.

It follows from \cite{KL}, Theorem 1.1, that $G$ fits into a short exact sequence
$$ 1 \rightarrow H \rightarrow G \rightarrow L \rightarrow 1$$
where $H$ is virtually $\Z$ and $L$ is a uniform lattice in the
isometry  group of $T$, and thus Gromov hyperbolic.

The short exact sequence above is obtained by constructing a
quasi-action  of $G$ on $T \times \R$ so that $L$ is the image and
$H$ is the kernel of this quasi-action.  It then follows from Lemma
\ref{virtual-center} that $H$ is the virtual center of $G$, and thus
characteristic under any group automorphism $\vp \in Aut(G)$. This
allows us to induce a surjective homomorphism $\overline{\vp}: L
\rightarrow L$.  Since $L$ is a non-elementary Gromov hyperbolic
group, it is Hopfian by \cite{S} and so Lemma \ref{lemma:hopfian}
implies that $\overline{\vp}\in Aut(L)$. It follows from \cite{LL,F}
that $R(\overline \vp )=\infty$, and from Lemma \ref{lemma:reid}
that $R(\vp) = \infty$ as well.

We note that this case could also be proven using the main result of
\cite{R},  which is a special case of Theorem 1.1 of \cite{KL}.

{\bf Case 3.} If $G$ is quasi-isometric to $\Gamma$, which is
quasi-isometric  to $BS(2,3)$, then $G$ itself is quasi-isometric to
$BS(2,3)$.  We quote the following theorem of Whyte which describes
groups quasi-isometric to $BS(2,3)$.

\begin{theorem}[\cite{W}, Theorem 5.1]\label{whyte}
Let $\Gamma$ be a finitely generated group.  Then $\Gamma$ is
quasi-isometric  to $BS(2,3)$ iff $\Gamma$ is the fundamental group
of a graph of virtual $\Z$'s which is neither commensurable to $F_n
\times \Z$ nor virtually solvable.
\end{theorem}

Thus our group $G$ is the fundamental group of a graph of groups all
of whose  vertex and edge groups are virtually infinite cyclic.  Let
$T$ be the tree on which $G$ acts with stabilizers which are
virtually infinite cyclic.  The elliptic elements of $G$ fall into
two classes: those with finite order, and those with infinite order.
Using the fact that all infinite order elliptic elements in $G$ are
commensurable, we can define the modular homomorphism $\Delta: G
\rightarrow \Q^*$ as follows.

Fix an infinite order elliptic element $\alpha$, and let $g$ be any
element of $G$.  Since $g \alpha g^{-1}$ will be an infinite order
elliptic element, and thus commensurable to $\alpha$, we see that
there is a relator of the form $g \alpha^p g^{-1} = \alpha^q$ for
some $p,q \in \Z - \{0\}$.  Define $\Delta(g) = \frac{p}{q}$, which
is well defined because all infinite order elliptic elements are
commensurable. Since automorphisms of $G$ preserve both finite or
infinite order and the type of element (elliptic or hyperbolic), we
see that $\Delta \circ \psi = \Delta$ for any $\psi \in Aut(G)$.

We say that a group $G$ is {\em unimodular} if the image of $\Delta$
is  contained in $\{ \pm 1\}$.  We now quote Levitt's proof from
\cite{L}, Proposition 2.7 that when $G$ is not unimodular, $R(\vp) =
\infty$.  Namely, the image of $\Delta$ is infinite, and
$\vp$-conjugate elements have the same modulus, so the result
follows.

We now suppose that $G$ is unimodular, and adapt the proofs of
\cite{L}, Propositions 2.5, 2.6 and 2.7.  Following Levitt, we note
that all central elements in $G$ are elliptic, and that the center
$Z(G)$ is contained in the kernel of the action of $G$ on $T$.  In
the case of generalized Baumslag-Solitar groups, one can conclude
that $Z(G)$ is trivial or infinite cyclic \cite{L}.  However, in our
case, we can only say that it is either virtually $\Z$ or the normal
closure of the torsion elements of $G$, since $Z(G)$ must be
contained in every vertex stabilizer, all of which are virtually
$\Z$.

If $\Delta(G)$ is trivial, we conclude following Levitt's argument
that $Z(G)$ is virtually $\Z$.  In this case, we note that $G/Z(G)$
acts on $T$ with finite stabilizers, and thus must be virtually
free.  We obtain the short exact sequence
\begin{equation}\label{unimodular}
1 \rightarrow Z(G) \rightarrow G \rightarrow F \rightarrow 1
\end{equation}
where $F$ is virtually free.  Since we have assumed that $T$ is not
a line, we know that $F$ contains a free group $F_n$ with $n \geq 2$
as a subgroup  of finite index. The pullback of \eqref{unimodular}
by the inclusion $F_n\hookrightarrow F$ is an extension $G'$ of
$Z(G)$ by $F_n$, namely
\begin{equation}\label{unimodular2}
1 \rightarrow Z(G) \rightarrow G' \rightarrow F_n \rightarrow 1
\end{equation}
Since $F_n$ is free, the sequence splits. When we view the group
$G'$ as a finite index subgroup of $G$, and recall that $Z(G)$ is
the center of $G$, it follows that $G'$ is a direct product, that
is, $G'=Z(G)\times F_n$. From this we conclude that $G$ is virtually
$Z(G)\times F_n$.  As $Z(G)$ is virtually infinite cyclic, we see
that $G$ is virtually $\mathbb Z \times F_n$. Following Theorem
\ref{whyte}, the case $\Delta(G)=1$ cannot occur.

We now consider the case when $Z(G)$ is the normal closure of the
torsion  elements of $G$, which corresponds to the case when
$\Delta(G) \subset \{\pm 1\}$.  Consider $K = Ker(\Delta)$ which has
index two in $G$.  We note that $K$ is also the fundamental group of
a graph of groups all of whose vertex and edge groups are virtually
$\Z$: $K$ acts on the same tree as $G$, and $[G:K] = 2$ implies that
all vertex and edge stabilizers  must again be virtually infinite
cyclic. Since $\Delta(K)$ is trivial, we conclude as in the case
$\Delta(G)=1$ that this case cannot occur. Hence, in case 3, $G$
cannot be unimodular and the proof is complete.
\end{proof}


\begin{thebibliography}{TW}

\bibitem[B]{B}
R. Brown et al (eds.), Handbook of topological fixed point theory, Springer, 2005.

\bibitem[EFW1]{EFW1}
A. Eskin, D. Fisher, and K. Whyte, Coarse differentiation of
quasi-isometries I: spaces not quasi-isometric to Cayley graphs,
preprint, 2007.

\bibitem[EFW2]{EFW2}
A. Eskin, D. Fisher, and K. Whyte, Coarse differentiation of
quasi-isometries II: Rigidity for Sol and Lamplighter groups,
preprint, 2007.

\bibitem[F]{F}
A.L. Fel'shtyn, The Reidemeister number of any automorphism of a
Gromov hyperbolic group is infinite, {\it Zap. Nauchn. Sem. POMI}
Vol. 279(2001), 229-241.

\bibitem[FG1]{FG}
A.L. Fel'shtyn and D. Gon\c calves,  Twisted conjugacy classes of
automorphisms of Baumslag-Solitar groups, {\it Algebra Discrete
Math.} No. 3 (2006), pp. 36-48.

\bibitem[FG2]{FG2}
A.L. Fel'shtyn and D. Gon\c calves, Twisted conjugacy classes in
mapping class groups, symplectic groups and braid groups, (Appendix:
Geometric group theory and $\Rinf$ property for mapping class
grouops, with F. Dahmani), preprint, 2007.

\bibitem[GW1]{GW1}
D. Gon\c calves and P. Wong, Twisted conjugacy classes in wreath
products,  {\it Internat. J. Alg. Comput.}, Vol. 16 No. 5 (2006),
pp. 875-886.

\bibitem[GW2]{GW2}
D. Gon\c calves and P. Wong, Twisted conjugacy in exponential growth
groups,  {\it Bull. London Math. Soc.}, Vol. 35 (2003), 261-268.

\bibitem[GW3]{GW3}
D. Gon\c calves and P. Wong, Twisted conjugacy classes in nilpotent
groups, {\em J. Reine Angew. Math.}, to appear.

\bibitem[J]{J}
B. Jiang, Lectures on Nielsen Fixed Point Theory, Contemp. Math. {\bf 14}, Amer. Math. Soc., Providence 1983.

\bibitem[KL]{KL}
B. Kleiner and B. Leeb, Groups quasi-isometric to symmetric spaces, {\it Comm. Anal. Geom.} 9 No. 2 (2001) pp. 239-260.

\bibitem[L]{L} G. Levitt, On the automorphism group of generalized Baumslag-Solitar groups,
{\it Geometry and Topology} 11 (2007), 473--515.

\bibitem[LL]{LL}
G. Levitt and M. Lustig, Most automorphisms of a hyperbolic group have simple dynamics, {\it Ann. Sci. Ecole Norm. Sup.} 33 (2000), 507-517.

\bibitem[MSW]{MSW}
L. Mosher, M. Sageev and K. Whyte, Quasi-actions on trees I, {\it Ann. of Math.}(2) 158 (2003), 115-164.

\bibitem[R]{R}
E. G. Rieffel, Groups Quasi-isometric to the Hyperbolic Plane Cross the Real Line. {\it Journal of the London Mathematical Society}, Vol. 64 no. 1(2001), pp. 44 - 60.

\bibitem[S]{S}
Z. Sela, Endomorphisms of hyperbolic groups I: The Hopf property, {\it Topology} 38 (1999), 301-321.

\bibitem[TWh]{TWh}
J. Taback and K. Whyte, Twisted conjugacy and group actions, preprint, 2005.

\bibitem[TWo1]{TW1}
J. Taback and P. Wong, Twisted conjugacy and quasi-isometry
invariance for  generalized solvable Baumslag-Solitar groups,
{\it Journal London Math. Soc. (2)} 75 (2007), 705-717.

\bibitem[TWo2]{TW2}
J. Taback and P. Wong, The geometry of twisted conjugacy classes in
wreath products, preprint, 2008.

\bibitem[W]{W}
K. Whyte,  The large scale geometry of the higher Baumslag-Solitar groups.  {\it Geom. Funct. Anal.}  11  ,  no. 6 (2001), pp.  1327--1343.

\bibitem[Wo]{Wo}
P. Wong. Reidemeister number, Hirsch rank, coincidences on polycyclic groups and solvmanifolds. {\it J. reine angew. Math.}, 524(2000), 185-204.

\bibitem[Wo2]{Wo2}
P. Wong, Fixed point theory for homogeneous spaces---a brief survey.  Handbook of topological fixed point theory,  265--283, Springer, Dordrecht, 2005.

\bibitem[Wo3]{Wo3}
P. Wong, Fixed point theory for homogeneous spaces. II. {\it Fund. Math.}  186  (2005),  no. 2, 161--175.

\end{thebibliography}
\end{document}